\documentclass[11pt]{article}

\usepackage[utf8]{inputenc}
\usepackage[T1]{fontenc}
\usepackage{amsmath}
\usepackage{amsfonts}
\usepackage{amssymb}
\usepackage{graphicx}

\usepackage{authblk}
\usepackage{mathrsfs}

\begin{document}

\title{\bf A New Overture to Classical Simple Type Theory, Ketonen-type Gentzen and Tableau Systems\footnote{ arXiv edition, The first version, v.1.0}}

\author[1]{\Large Tadayoshi Miwa \thanks{miwa.tadayoshi@mail.u-tokyo.ac.jp}}
\author[2]{\Large Takao Inou\'{e} \thanks{Corresponding Author: takaoapple@gmail.com}}

\affil[1]{\normalsize Library, The University of Tokyo, Tokyo, Japan}
\affil[2]{Faculty of Informatics, Yamato University, Osaka, Japan\footnote{inoue.takao@yamato-u.ac.jp}}


\date{January 15, 2026}


\maketitle

\begin{abstract}
In this paper, we introduce a Ketonen-type Gentzen-style classical simple type theory $\bf KCT$. Also the tableau system $\bf KCTT$ corresponding to $\bf KCT$ is introduced. Further inference-preserving Gentzen system $\bf KCT_h$ (equivalent to $\bf KCT$) and tableau system $\bf KCTT_h$ (equivalent to $\bf KCTT$) is introduced. We introduce the notion of Hintikka sequents for $\bf KCTT_h$.The completeness theorem and Takahashi-Prawitz's theorem are proved  for $\bf KCTT_h$. 
\medskip

\noindent MSC2020: 03B38, 03F03
\medskip

\noindent \small \it Keywords: \rm Ketonen-type classical simple type theory, Ketonen-type Gentzen systems, Ketonen-type Tableau methods, inference-preserving system, Hintikka sequent, partial valuation, Prawitz's extension of partial valuation, completeness theorem,  Takahashi-Prawitz's theorem.

\end{abstract}
	
\tableofcontents

 \normalsize
	
\newtheorem{theorem}{Theorem}[section]
\newtheorem{proposition}{Proposition}[section]
\newtheorem{definition}{Definition}[section]
\newtheorem{corollary}{Corollary}[section]
\newtheorem{lemma}{Lemma}[section]
\newtheorem{convention}{convention}[section]
\newtheorem{remark}{Remark}[section]


\section{Introduction}
\subsection{General introduction: Proposed systems and the results}

This paper is organized as follows. After this general introduction, we first explain the significance of simple type theory along historical line with technical details in \S 1.2. In \S 2 and \S 3, we give formal preliminaries for our systems.  We introduce a Ketonen-type Gentzen-style classical simple type theory $\bf KCT$ in \S 4.  Also in \S 5 the tableau system $\bf KCTT$ corresponding to $\bf KCT$ is introduced. Further inference-preserving Gentzen system $\bf KCT_h$ (equivaltent to $\bf KCT$) in \S 6 and tableau system $\bf KCTT_h$ (equivalent to $\bf KCTT$) is introduced in \S 7. A number of definitions are introduded for the proofs of our main theorems in \S 8: among others, reduction chains and partical valuations for formulas and sequents. Some lemmata are presented in \S 9. From \S 9 to the last, our interest moves to the tableau system  $\bf KCTT_h$. We introduce the notion of Hintikka sequents for $\bf KCTT_h$ in \S 10. We mention some deducible sequents and permissible inferences for $\bf KCTT_h$. The semanics is prepared in \S 12. The completeness theorem and Takahashi-Prawitz's theorem are proved  for $\bf KCTT_h$ (cf. Takahashi \cite{takahashi1967} and Prawitz \cite{prawitz1968}) in \S 13. We conclude this paper, that is, a new overture to simple type theory with a short remark. 

We follow the basic notations and notions of Sch\"{u}tte \cite{schutte1977} for $\bf KCT$, $\bf KCTT$, $\bf KCT_h$ and $\bf KCTT_h$. 

\paragraph{Main contributions.}
The main contributions of this paper are as follows:
\begin{itemize}
  \item We introduce a Ketonen-type Gentzen-style simple type tableau system formulated as an \emph{inference-preserving system} $\bf KCTT_h$, rather than merely as a decision procedure.
  \item We define and systematically use \emph{Hintikka sequents}  (\S 10)  as structural objects mediating between tableau rules and sequent-style reasoning.
  \item We show that the proposed system provides a uniform proof-theoretic framework for classical simple type theory, thereby clarifying the relationship between tableau calculi and Gentzen systems.
  \item We show that the presentation of simple type theory in Sch\"{u}tte's Proof Theory \cite[Chapter IV]{schutte1977} can be naturally understood as a Ketonen-type tableau system, and we make this correspondence explicit.
\end{itemize}

Using the Ketonen-type framework, 
we establish completeness results for the proposed systems by means of inference-preserving cut elimination, introducing the notion of 
Hintikka sequents. This paper can also be considered a reconstruction of Schutte's proof theory approach using the Ketonen tableau method.

\subsection{Why Type Theory?---Russell's Paradox, Syntactic Stratification, and Foundations}
\label{subsec:why-type-theory}

Type theory was introduced as a principled response to paradoxes arising from
\emph{unrestricted comprehension} principles in naive set-theoretic and
propositional reasoning, most famously \textbf{Russell's paradox} \cite{Russell1908}.
Consider the comprehension
\[
  R = \{\, x \mid x \notin x \,\}.
\]
If $R \in R$, then by definition $R \notin R$; conversely, if $R \notin R$, then
$R \in R$, yielding a contradiction.

The core idea behind type theory is to \textbf{stratify expressions into levels
(types)} so that illegitimate forms of self-application become ill-typed and hence
ill-formed. To restore consistency at the foundational level, Russell introduced
the concept of \textbf{stratification}, motivated by the \emph{vicious circle
principle}, ensuring that predicates have a restricted range of significance
\cite{Russell1908}.
In his \emph{ramified} theory of types, developed in the tradition of
\emph{Principia Mathematica} \cite{WhiteheadRussell1910}, logical entities are
arranged in a hierarchy of both types and orders, where a function (or predicate)
may take arguments only of an appropriate lower type and order.
This stratification enforces a strict syntactic discipline that excludes
expressions such as $x(x)$ or, in a set-theoretic guise, $x \in x$ as ill-formed,
rather than merely false, thereby structurally preventing the formulation of
self-referential paradoxes.

Church's \textbf{simple type theory} \cite{church1940} distills this stratification
idea into a streamlined and technically manageable framework based on function
types.
In Church's formulation, which provides the typed logical background for the system
studied in this paper, terms are assigned specific types, 
and logical operations are defined relative to this type hierarchy.
Consequently, a term of a functional type may be applied only to arguments of the
appropriate input type, and quantification is explicitly type-indexed.
Consistency is thus safeguarded at the syntactic level, in the sense that the
grammatical formation rules block the impredicative constructions required to
formulate Russell-style paradoxes.

The calculi investigated in this paper, including Gentzen-style systems and their
Ketonen-type variants, are formulated over such a typed language.
Our proof-theoretic results, such as completeness and cut elimination, are therefore
established within a framework in which Russell-style self-reference is excluded at
the level of well-formedness.
In this way, the operational virtues of sequent calculi and inference-preserving
systems are aligned with the foundational motivation that originally gave rise to
type-theoretic discipline.

\section{Preparations for Our Formal Systems}

\subsection{The Formal Language}

We mainly follow Sch\"{u}tte's  terminology and definitions for simple type theory (Refer to Sch\"{u}tte \cite{schutte1977}).

As usual in \cite{schutte1977}, by an \bf $n$-place nominal form \rm ($n \geq 1$) we mean a finite string of 
symbols in which besides the primitive symbols of our formal language (sentential variables, connectives and 
round brackets) the \bf nominal symbols \rm $*_1, \dots , *_n$ occur just once. If $\mathscr{E}$ is an $n$-simple nominal form ($n \geq 1$) and $r_1, \dots , r_n$ are non-empty finite strings of symbols then $\mathscr{E}[r_1, \dots , r_n]$ denotes the result of replacing the nominal symbols  $*_1, \dots , *_n$ in  $\mathscr{E}$ by $r_1, \dots , r_n$, respectively.

\begin{definition}
Inductive definition of the \textbf{types} is as follows:

2.1.1. \(0\) and \(1\) are types. 

2.1.2. If \(\tau_1, \ldots, \tau_n\) \((n \geq 1)\) are types, then \((\tau_1, \ldots, \tau_n)\) is also a type. 
\end{definition}

We use \(\sigma, \sigma_i, \tau_i\) as \textbf{syntactic symbols} for types and \(\tau\) for sequences \(\tau_1, \ldots, \tau_n\) \((n \geq 1)\) of types.

\begin{definition}
\textbf{Primitive Symbols} of our System are as follows:

2.2.1 Denumerably infinitely many free and bound variables of each type. 

2.2.2 A non-empty set of object symbols. 

2.2.3 Certain function letters with specified numbers \(\geq 1\) of arguments. 

2.2.4 The symbols \(\supset\), \(\forall\) and \(\lambda\). 

2.2.5 Round brackets and comma.
\medskip 

For \textbf{syntactic symbols} for variables we use
\begin{itemize}
  \item $a^\sigma, a_i^{\tau_i} \quad \text{for free variables of types } \sigma, \tau_i,$
  \item $x^\sigma, x_i^{\tau_i} \quad \text{for bound variables of types } \sigma, \tau_i,$
  \item $a^{\tau} \quad \text{ for sequences } a_1^{\tau_1} \dots a_n^{\tau_n} (n \geq 1)$ of free variables,
  \item $x^{\tau} \quad \text{ for sequences } x_1^{\tau_1} \dots x_n^{\tau_n} (n \geq 1)$ of pairwise distinct bound variables.
\end{itemize}
\end{definition}

\begin{definition}
Inductive definition of \textbf{terms} in the system and their types is the following:

2.3.1 Every free variable of type \(\sigma\) is a \textbf{term} of type \(\sigma\). 

2.3.2 Every object symbol is a term of type 0. 

2.3.3 If \(\varphi\) is an \(n\)-place function letter \((n \geq 1)\) and \(t_1^0, \ldots, t_n^0\) are terms of type 0, then \(\varphi(t_1^0, \ldots, t_n^0)\) is a term of type \(0\). 

2.3.4 If \(t_1^{\tau_1}, \ldots, t_n^{\tau_n} \ (n \geq 1)\) and \(t_0^{({\tau_1}, \ldots, {\tau_n})}\) are terms of types \(\tau_1, \ldots, \tau_n\) and \((\tau_1, \ldots, \tau_n)\), 

\noindent then
  \(t_0^{({\tau_1}, \ldots, {\tau_n})} (t_1^{\tau_1}, \ldots, t_n^{\tau_n})\) is a term of type 1. 
  
 $ ($We sometimes denote this by \(t_0^{(\tau)}(t^{\tau})\) for brevity where \(\tau = \tau_1, \ldots, \tau_n\) and $t^{\tau} = (t_1^{\tau_1}, \ldots, t_n^{\tau_n}).)$
  
2.3.5 If \(A\) and \(B\) are terms of type 1 then \((A \supset B)\) is also a term of type 1.

2.3.6 If \(\mathscr{F}[a^\sigma]\) is a term of type 1 and \(x^\sigma\) is a bound variable which does not occur in the nominal form \(\mathscr{F}\) then \(\forall x^{\sigma} \mathscr{F}[x^{\sigma}]\) is also a term of type 1. 

2.3.7 If \(\mathscr{A}[a_1^{\tau_1}, \ldots, a_n^{\tau_n}]\ (n \geq 1)\) is a term of type 1 and \(x_1^{\tau_1}, \ldots, x_n^{\tau_n}\) are pairwise distinct bound variables which do not occur in the \(n\)-place nominal form \(\mathscr{A}\) then \(\lambda x_1^{\tau_1}, \ldots, x_n^{\tau_n} \mathscr{A}[x_1^{\tau_1}, \ldots, x_n^{\tau_n}]\) is a term of type $(\tau_1, \ldots, \tau_n)$. $($We sometimes denote this by $\lambda x^\tau \mathscr{A}[x^\tau]$ for brevity where $x^\tau = x_1^{\tau_1}, \ldots, x_n^{\tau_n}$.$)$
\end{definition}

\begin{definition}
The terms of type $1$ are called \textbf{formulas}. Free variables of type $1$ and formulas of the form $a_0^{(\tau_1, \ldots, \tau_n)} (t_1^{\tau_1}, \ldots, t_n^{\tau_n})$ where $a_0^{(\tau_1, \ldots, \tau_n)}$ is a free variable, are called \textbf{atomic formulas}.
\end{definition}

For \textbf{syntactic symbols} for the rest we use
\begin{itemize}
  \item $t^{\sigma}$, $t_i^{\tau_i}$ \quad for terms of types $\sigma$, $\tau_i$,
  \item $t^{\tau}$ \quad for sequences $t_1^{\tau_1}, \ldots, t_n^{\tau_n}$ ($n \geq 1$) of terms,
  \item $A$, $A_i$, $B$, $B_i$, $C$, $C_i$, $F$, $F_i$, $G$, $G_i$ \quad for formulas,
  \item $P$, $P_i$ \quad for atomic formulas,
  \item $\mathscr{F}$, $\mathscr{F}_i$ \quad for 1-place nominal forms such that $\forall x^{\sigma} \mathscr{F} [x^{\sigma}]$, $\forall x^{\sigma} \mathscr{F}_i [x^{\sigma}]$
  are formulas,
  \item $\mathscr{A}$, $\mathscr{A}_i$ \quad for $n$-place nominal forms ($n \geq 1$) such that $\lambda x^\tau \mathscr{A} [x^\tau], \lambda x^\tau \mathscr{A}_i [x^\tau]$ are terms of type $(\tau)$ where $\tau = \tau_1, \ldots, \tau_n$.
\end{itemize}

We identify terms which only differ in the choice of the bound variables occurring in them. When we write $\mathscr{F} [t^\sigma]$ or $\mathscr{A} [t^\tau]$ we assume that the bound variables which occur in $t^\sigma$ and $t^\tau$ are chosen so that they do not also occur in the nominal form $\mathscr{F}$ or $\mathscr{A}$ (respectively).

\begin{definition} 

We set $ \bot := \forall x^1x^1$. For any formulas $A, B$, we set  

2.5.1.  $ \neg A := A \supset \bot,$
 
2.5.2. $A \vee B := \neg A \supset B,$

2.5.3. $A \wedge B := \neg (A \supset \neg B),$

2.5.4. $\exists x A := \neg \forall x \neg A.$
\end{definition}

\begin{definition} 
We define a \textbf{sequent} as $\Gamma \rightarrow \Delta$, where $\Gamma, \Delta$ are sequents of formulas,  and $\Gamma, \Delta$ may be empty sequents. Let $\Gamma \rightarrow \Delta$ be a sequent where $\Gamma, \Delta$ are sequents of formulas. Every formula contained in $\Gamma$ is called  to be \bf A-formula \rm and every formula contained in $\Delta$ is called to be \bf S-formula\rm .
\end{definition}

The intuitive meaning of a sequent is as follows. 

IM1. $A_1, A_2, \dots, A_n \rightarrow B_1, B_2, \dots , B_m$ 

$\qquad \qquad (n \geq  1, m \geq 1) \quad$ is 

$(A_1 \wedge A_2 \wedge \dots \wedge A_n) \supset (B_1 \vee B_2 \vee  \dots \vee B_m).$
\smallskip

IM2. $A_1, A_2, \dots , A_n \rightarrow $ $(n \geq  1) \quad$ is 

$\qquad \neg (A_1 \wedge A_2 \wedge \dots \wedge A_n).$
\smallskip

IM3. $\rightarrow B_1, B_2, \dots , B_m$ $(m \geq 1) \quad$ is 

$\qquad B_1 \vee B_2 \vee  \dots \vee B_m.$
\smallskip

IM4. $\rightarrow $ is $\bot$.

\begin{definition} 
By \textbf{the corresponding formula} of a sequent $S$ we mean the formula $F_S$ which is the intuitive meaning of $S$ interpreted by IM1-IM4.
\end{definition} 

Following Sch\"{u}tte \cite{schutte1977}, we identify terms which only differ in the choice of the bound variables occurring in them. When we write $\mathscr{F}[t^\sigma]$ or $\mathscr{A}[t^\tau]$ we suppose that the bound variables which occur in $t^\sigma$ and $t^\tau$ are chosen so that they do not also occur in the nominal form $\mathscr{F}$ or  $\mathscr{A}$, respectively.

\begin{proposition} $ $

1.  If \( \forall x^{\sigma} \mathscr{F}[x^\sigma] \) is a formula and $t^{\sigma}$ is a term of type $\sigma$, then \( \mathscr{F}[t^\sigma] \) is a formula. 

2.  If $\lambda x^\tau \mathscr{A} [x^\tau]$ is a term of type \( (\tau) = (\tau_1, \ldots, \tau_n) \) and if $t^\tau = t_1^{\tau_1}, \ldots, t_n^{\tau_n}$ is a sequence of terms of types \( \tau_1, \ldots, \tau_n \) then \( \mathscr{A} [t^\tau] \) is a formula. 
\end{proposition}

\subsection{Chains of Subterms}

In this subsection, we need to give formulas a rank for our proofs by induction. We define this rank with the aid of chains of subterms.

\begin{definition} 
We call terms of type 0 and free variables \textbf{prime terms}.
\end{definition} 

\begin{definition} 
We define the \textbf{subterms} of a term as follows:

2.9.1 A prime term has no subterms.

2.9.2 An atomic formula $a_0^{(\tau_1, \ldots, \tau_n)} (t_1^{\tau_1}, \ldots, t_n^{\tau_n})$ has $t_1^{\tau_1}, \ldots, t_n^{\tau_n}$ as subterms.

2.9.3 A formula \( \lambda x^\tau \mathscr{A} [x^\tau] (t^\tau) \) has $\mathscr{A} [t^\tau]$ as subterm.
  
2.9.4 A formula $(A \supset B)$ has \( A \) and \( B \) as subterms.
  
2.9.5 A formula $\forall x^\sigma \mathscr{F} [x^\sigma]$ has as subterms every formula $\mathscr{F} [a^\sigma]$ (for every free variable \( a^\sigma \) of type \( \sigma \)).
  
2.9.6 A term $\lambda x^\tau \mathscr{A} [x^\tau]$ has as subterms every formula $\mathscr{A} [a^\tau]$ (for every sequence $a^\tau = a_1^{\tau_1}, \ldots, a_n^{\tau_n}$ of free variables of types $\tau_1, \ldots, \tau_n$).

\end{definition} 

\begin{definition}
By a \textbf{subterm chain} of a term $t^\sigma$ we mean a sequence
$$
t_0^{\tau_0}, t_1^{\tau_1}, t_2^{\tau_2}, \ldots
$$
of terms formed as follows:

2.10.1 The initial term $t_0^{\tau_0}$ of the subterm chain is the term $t^\sigma$.

2.10.2 If a term $t_i^{\tau_i}$ in the subterm chain is not a prime term, then the subterm chain contains as immediate successor of $t_i^{\tau_i}$ a subterm $t_{i+1}^{\tau_{i+1}}$ of $t_i^{\tau_i}$ .
  
2.10.3 If a term $t_i^{\tau_i}$ in a subterm chain is a prime term then it is the last term in the subterm chain. We then say that the subterm chain has \textbf{length} $i$.

\end{definition} 

\begin{definition} 
A term $t^\sigma$ is said to be \textbf{regular} if all subterm chains of $t^\sigma$ have a maximal finite length $m$. We define this number $m$ to be the \textbf{rank} $Rt^\sigma$ of a regular term $t^\sigma$.
\end{definition} 

In the following we suppose $\nu$ is a $1$-place nominal form such that $\nu[a^\sigma]$ is a term.

\begin{lemma}\textup{(\cite{schutte1977})}

If $\nu[a_1^{\sigma}]$ is a regular term, then $\nu[a_2^{\sigma}]$ is a regular term with the same rank as $\nu[a_1^{\sigma}]$.
\end{lemma}

\begin{lemma}\textup{(\cite{schutte1977})} 
If every subterm of a term $t^{\sigma}$ is regular, then $t^{\sigma}$ is regular.
\end{lemma} 

\begin{definition}
We define the \textbf{height} of a type $\sigma$ to be the number of round brackets which occur in $\sigma$.
\end{definition}

\begin{lemma}\textup{(\cite{schutte1977})}

The following conditions are sufficient for a term $\nu[t^\sigma]$ to be regular:

\textup{(1)} $\nu[a^\sigma]$ and $t^\sigma$ are regular terms.

\textup{(2)} For every type $\tau_i$ of height less than $\sigma$ we have: If $\mathscr{W}[a_i^{\tau_i}]$ and $t_i^{\tau_i}$ are regular terms, then $\mathscr{W}[t_i^{\tau_i}]$ is a regular term.

\textup{(3)} For every regular term $\nu_i[\alpha^\sigma]$ of smaller rank than $s[a^\sigma]$, $\nu_i [t^\sigma]$ is a regular term.
\end{lemma} 

\begin{lemma}\textup{(\cite{schutte1977})} If $\nu[a^\sigma]$ and $t^\sigma$ are regular terms, then $\nu[t^\sigma]$ is also a regular term.
\end{lemma}
\begin{definition}
We define the \bf length \rm of a term $t^\sigma$ to be the number of primitive symbols which occur in $t^\sigma$
\end{definition}
\begin{theorem}\textup{(\cite{schutte1977})} Every term $t^\sigma$ is regular.
\end{theorem}

\begin{corollary} \textup{(\cite{schutte1977})} Every formula has a computable finite rank.
\end{corollary}

\begin{theorem}\textup{(\cite{schutte1977})}

The following hold for ranks of formulas:

\textup{(1)} $R A < R (A \supset B) \quad \text{and} \quad R B < R (A \supset B)$.

\textup{(2)} $R \mathscr{F} [a^{\sigma}] < R \forall x^{\sigma} \mathscr{F} [x^{\sigma}]$.

\textup{(3)} $R \mathscr{A}[ t^\tau] < R \lambda x^\tau \mathscr{A} [ x^\tau ] ( t^\tau )$.
\end{theorem}

\section{Basic preliminaries on inferences}
\paragraph{(1) Basic inferences}
Certain configurations of the form $$A_1, \ldots, A_n \vdash B $$ where $A_1, \ldots, A_n, B$ are formulas and in general $n=1$ or $n=2$, are defined to be  \textbf{basic inferences}. A formula on the left of the symbol $\vdash$ is said to be a \textbf{premise} of the given basic inference. The formula on the right of the symbol $\vdash$  is said to be the \textbf{conclusion} of the given basic inference.

\paragraph{(2) deducibility}
Deducibility of formulas in a formal system is determined by the axioms and
basic inferences according to the following inductive definition.
\begin{enumerate}
\item Every axiom is \emph{deducible}.
\item If every premise of a basic inference is \emph{deducible}, then the conclusion of the inference is \emph{deducible}.
\end{enumerate}
A proof (demonstration) of the deducibility of a formula \(A\) (according to 1 and 2 above) is said to be a \textbf{deduction} of \(A\).

\begin{definition}[Order of Deducibility]
Sometimes we need to specify the deducibility of a formula more closely; so to
the formula we assign a natural number, which serves as the \emph{order}, by the following
inductive definition:

DO 1. Every axiom is \emph{deducible} with \emph{order} $0$.

DO 2. If the premises $A_1,\ldots,A_n$ $(n\ge 1)$ of a basic inference are deducible with
  orders $m_1,\ldots,m_n$, respectively, then the conclusion of the inference is deducible
  with order
  $$\max(m_1,\ldots,m_n)+1$$
A formula may be deducible with many different orders. However, given that a
formula is deducible we may always assume that it is deducible with some fixed
order.
\end{definition}

\paragraph{(3) reducibility}
As it were, the reduction is the converse of deducibility. Given a formula $A$ to be provable, $A$ is reduced to axioms by application of reduction rules.

In \S 7 we present such a tableau system as our main method. In \S 8, we fix the notions of reducible sequents and
reduction chains. In this paper, our main theoretical interest lies in reducibility.

\paragraph{(4) The turnstile \(\vdash\) in this paper (a caution).}
Throughout the paper we use the form \(\vdash_{\mathsf{X}}\) as a \emph{derivability relation on sequents} (in Gentzen or tableau), not primarily on formulas.
Specifically, for any system \(\mathsf{X}\in\{\mathsf{KCT},\mathsf{KCT}_{h},\mathsf{KCTT},\mathsf{KCTT}_{h}\}\) and any sequent \(S\), we write $\vdash_{\mathsf{X}}\,S$ to mean that there exists a finite derivation (proof figure or closed tableau, as appropriate) whose \emph{root} is \(S\).
When a single formula \(A\) is said to be `provable' in \(\mathsf{X}\), this is a shorthand for the relation $\vdash_{\mathsf{X}} A$, 
i.e., $\vdash_{\mathsf{X}} \enspace \rightarrow A$ in the sequent form. One may omit the $X$ if one knows the context well.

\paragraph{(6) permissible inference}
A configuration $A_1, \dots, A_n \vdash B$ 
where \(A_1, \dots, A_n\) and \(B\) are formulas (or sequents) \((n \ge 1)\) is said to be a
\emph{permissible inference} if the following holds: if \(A_1, \dots, A_n\) are
deducible, then \(B\) is also deducible.

\paragraph{(7) weak inference}
A \emph{weak inference} is a permissible inference \(A \vdash B\) such that:
if the formula (or sequent) \(A\) is deducible with order \(m\), then the formula (or sequent, resp.) \(B\)
is deducible with order \(\le m\).

\section{Ketonen-type Gentzen classical simple type theory $\bf KCT$}

Here we introduce a Ketonen-type Gentzen classical simple type theory $\bf KCT$ as follows.
\medskip

Axioms: $\Gamma_1, P, \Gamma_2 \enspace \rightarrow \enspace \Delta_1, P, \Delta_2,$
\medskip

\noindent where $P$ is atomic and all A-formulas and S-formula of $\Gamma_1, \Gamma_2, \Delta_1, \Delta_2$ are atomic and $\Gamma_1, \Gamma_2, \Delta_1, \Delta_2$ may be empty, respectively.
\medskip

Rules of inference: 
\bigskip

$\displaystyle{\frac{\enspace \Gamma_1, A, \Gamma_2 \enspace \rightarrow \enspace \Delta_1, B, \Delta_2 \enspace}{\enspace \Gamma_1, \Gamma_2 \enspace \rightarrow \enspace \Delta_1, A \supset B, \Delta_2}}$ $ (\rightarrow \supset)$

\bigskip

$\displaystyle{\frac{\Gamma_1, A \supset \forall x^1x^1, \Gamma_2 \rightarrow \Delta \enspace \enspace \Gamma_1, B, \Gamma_2 \rightarrow \Delta \enspace }{\enspace \Gamma_1, A \supset B, \Gamma_2 \enspace \rightarrow \enspace \Delta}}$ $(\supset \rightarrow)$,

\medskip

\noindent if $B$ is not a formula $\forall x^1x^1$.

\bigskip

\indent $\displaystyle{\frac{\enspace \Gamma \enspace \rightarrow \enspace \Delta_1, \mathscr{F}[a^\sigma], \Delta_2 \enspace}{\enspace \Gamma \enspace \rightarrow \enspace \Delta_1, \forall x^\sigma \mathscr{F}[x^\sigma], \Delta_2 \enspace }}$ $ (\rightarrow \forall^\sigma)$,

\bigskip

\noindent provided that the variable $a^\sigma$ does not occur in the lower sequent.

\bigskip

$\displaystyle{\frac{\enspace \Gamma_1, \mathscr{F}[t^\sigma], \forall x^\sigma \mathscr{F}[x^\sigma], \Gamma_2 \enspace \rightarrow \enspace \Delta \enspace}{\enspace \Gamma_1,\forall x^\sigma \mathscr{F}[x^\sigma], \Gamma_2 \enspace \rightarrow \enspace \Delta \enspace }}$ $ (\forall^\sigma \rightarrow)$,
\bigskip

\noindent where $t^\sigma$ is a term.

\bigskip

$\displaystyle{\frac{\enspace \Gamma \enspace \rightarrow \enspace \Delta_1, \mathscr{A}[t^\tau], \Delta_2 \enspace }{\enspace \Gamma \enspace \rightarrow \enspace \Delta_1, \lambda x^\tau \mathscr{A}[x^\tau](t^\tau), \Delta_2 \enspace}}$ $(\rightarrow \lambda^\sigma)$

\bigskip

$\displaystyle{\frac{\enspace \Gamma_1, \mathscr{A}[t^\tau], \Gamma_2 \enspace \rightarrow \enspace \Delta \enspace}{\enspace \Gamma_1, \lambda x^\tau \mathscr{A}[x^\tau](t^\tau), \Gamma_2 \enspace \rightarrow \enspace \Delta \enspace }}$ $(\lambda^\sigma \rightarrow)$

\bigskip

The Ketonen-type Gentzen system to deduce contradictions for classical predicate logic was introduced in Inou\'{e} in \cite{Inoue1989}. The Ketonen-type Gentzen system for classical (resp. intuitionistic) predicate logic was introduced in Inou\'{e} in \cite{Inoue2024} (resp. in \cite{InoueK1}). 

\bigskip

In the system \(\mathsf{KCT}\), a \emph{proof figure} (derivation) is a finite tree whose leaves are instances of the axiom schema stated above and whose internal nodes are obtained by a single application of one of the inference rules given there. A sequent \(S\) is said to be provable in \(\mathsf{KCT}\) (written \(\vdash_{\mathsf{KCT}} S\)) iff there exists such
a derivation whose root is \(S\). In particular, a formula \(A\) is provable (written \(\vdash_{\mathsf{KCT}} A\)) iff there is a derivation whose end-sequent is the empty-antecedent sequent \(\,\to A\). Thus, whenever a proof figure that starts from axioms reaches a node of the form \(\to A\), the formula \(A\) has been proved in \(\mathsf{KCT}\).

\section{Ketonen-type tableau system for classical simple type theory $\bf KCTT$}

Here we introduce a Ketonen-type tableau system for classical simple type theory $\bf KCTT$ as follows.
\medskip

Axioms: $\Gamma_1, P, \Gamma_2 \enspace \rightarrow \enspace \Delta_1, P, \Delta_2,$
\medskip

\noindent where $P$ is atomic and all A-formulas and S-formula of $\Gamma_1, \Gamma_2, \Delta_1, \Delta_2$ are atomic and $\Gamma_1, \Gamma_2, \Delta_1, \Delta_2$ may be empty, respectively.
\medskip

Rules of inference:
\bigskip

$\displaystyle{\frac{\enspace \Gamma_1, \Gamma_2 \enspace \rightarrow \enspace \Delta_1, A \supset B, \Delta_2 \enspace }{\enspace \Gamma_1, A, \Gamma_2 \enspace \rightarrow \enspace \Delta_1, B, \Delta_2 \enspace}}$ $(\rightarrow \supset)$

\bigskip


$\displaystyle{\frac{\Gamma_1, A \supset B, \Gamma_2 \rightarrow \Delta \enspace }{\enspace \Gamma_1, A \supset \forall x^1x^1, \Gamma_2 \enspace \rightarrow \enspace \Delta \enspace \enspace \Gamma_1, B, \Gamma_2 \rightarrow \Delta}}$ $(\supset \rightarrow)$,

\medskip

\noindent if $B$ is not a formula $\forall x^1x^1$.


\bigskip

$\displaystyle{\frac{\enspace \Gamma \enspace \rightarrow \enspace \Delta_1, \forall x^\sigma \mathscr{F}[x^\sigma], \Delta_2 \enspace }{\enspace \Gamma \enspace \rightarrow \enspace \Delta_1, \mathscr{F}[a^\sigma], \Delta_2 \enspace}}$ $(\rightarrow \forall^\sigma)$,

\bigskip

\noindent provided that the variable $a^\sigma$ does not occur in the upper sequent.

\bigskip


$\displaystyle{\frac{\enspace \Gamma_1,\forall x^\sigma \mathscr{F}[x^\sigma], \Gamma_2 \enspace \rightarrow \enspace \Delta \enspace }{\enspace \Gamma_1, \mathscr{F}[t^\sigma], \forall x^\sigma \mathscr{F}[x^\sigma], \Gamma_2 \enspace \rightarrow \enspace \Delta \enspace}}$ $(\forall^\sigma \rightarrow)$,
\bigskip

\noindent where $t^\sigma$ is a term.

\bigskip

$\displaystyle{\frac{\enspace \Gamma \enspace \rightarrow \enspace \Delta_1, \lambda x^\tau \mathscr{A}[x^\tau](t^\tau), \Delta_2 \enspace}{\enspace \Gamma \enspace \rightarrow \enspace \Delta_1, \mathscr{A}[t^\tau], \Delta_2 \enspace }}$ $(\rightarrow \lambda^\sigma)$

\bigskip

$\displaystyle{\frac{\enspace \Gamma_1, \lambda x^\tau \mathscr{A}[x^\tau](t^\tau), \Gamma_2 \enspace \rightarrow \enspace \Delta \enspace }{\enspace \Gamma_1, \mathscr{A}[t^\tau], \Gamma_2 \enspace \rightarrow \enspace \Delta \enspace}}$ $(\lambda^\sigma \rightarrow)$

\bigskip

The Ketonen-type tableau system for classical predicate logic was introduced in Inou\'{e} in \cite{Inoue2024}. 

\bigskip

In the system \(\mathsf{KCTT}\), a \emph{tableau} for a sequent \(S\) is a finite tree of sequents whose root is \(S\), whose internal nodes are obtained by a single application of one of the tableau rules given there, and whose leaves are the resulting end-sequents. A branch is said to be \emph{closed} if its leaf is an instance of the axiom schema stated there; otherwise the branch is \emph{open}. A tableau is \emph{closed} if all of its branches are closed, and \emph{open} if it has at least one open branch. A tableau is called \emph{complete} (in the tableau sense) if, along every open branch, all applicable tableau rules have been exhaustively and fairly applied.

There is a well-known correspondence between analytic tableaux and sequent calculi: tableaux can be viewed as sequent systems upside-down. Hence a closed, complete tableau corresponds to a (cut-free) Gentzen-style proof figure, and conversely.

\section{Ketonen-type inference-preserving Gentzen classical simple type theory $\bf KCT_h$}
 
Here we introduce a Ketonen-type Gentzen classical simple type theory $\bf KCT_h$ as follows.
\medskip

Axioms: $\Gamma_1, P, \Gamma_2 \enspace \rightarrow \enspace \Delta_1, P, \Delta_2,$
\medskip

\noindent where $P$ is atomic and all A-formulas and S-formula of $\Gamma_1, \Gamma_2, \Delta_1, \Delta_2$ are atomic and $\Gamma_1, \Gamma_2, \Delta_1, \Delta_2$ may be empty, respectively.
\medskip

Rules of inference:
\bigskip

$\displaystyle{\frac{\enspace \Gamma_1, A, \Gamma_2 \enspace \rightarrow \enspace \Delta_1, B, A \supset B, \Delta_2 \enspace}{\enspace \Gamma_1, \Gamma_2 \enspace \rightarrow \enspace \Delta_1, A \supset B, \Delta_2}}$ $ (\rightarrow \supset)$

\bigskip

$\displaystyle{\frac{\Gamma_1, A \supset \forall x^1x^1, A \supset B, \Gamma_2 \rightarrow \Delta \enspace \enspace \Gamma_1, B,  A \supset B, \Gamma_2 \rightarrow \Delta \enspace }{\enspace \Gamma_1, A \supset B, \Gamma_2 \enspace \rightarrow \enspace \Delta \enspace }}$ 
$(\supset \rightarrow)$,

\medskip

\noindent if $B$ is not a formula $\forall x^1x^1$.

\bigskip

$\displaystyle{\frac{\enspace \Gamma \enspace \rightarrow \enspace \Delta_1, \mathscr{F}[a^\sigma],  \forall x^\sigma \mathscr{F}[x^\sigma], \Delta_2 \enspace}{\enspace \Gamma \enspace \rightarrow \enspace \Delta_1, \forall x^\sigma \mathscr{F}[x^\sigma], \Delta_2 \enspace }}$ $ (\rightarrow \forall^\sigma)$,

\bigskip

\noindent provided that the variable $a^\sigma$ does not occur in the upper sequent.

\bigskip

$\displaystyle{\frac{\enspace \Gamma_1, \mathscr{F}[t^\sigma], \forall x^\sigma \mathscr{F}[x^\sigma], \Gamma_2 \enspace \rightarrow \enspace \Delta \enspace}{\enspace \Gamma_1,\forall x^\sigma \mathscr{F}[x^\sigma], \Gamma_2 \enspace \rightarrow \enspace \Delta \enspace }}$ $ (\forall^\sigma \rightarrow)$,
\bigskip

\noindent where $t^\sigma$ is a term.

\bigskip

$\displaystyle{\frac{\enspace \Gamma \enspace \rightarrow \enspace \Delta_1, \mathscr{A}[t^\tau], \lambda x^\tau \mathscr{A}[x^\tau](t^\tau), \Delta_2 \enspace }{\enspace \Gamma \enspace \rightarrow \enspace \Delta_1, \lambda x^\tau \mathscr{A}[x^\tau](t^\tau), \Delta_2 \enspace}}$ $(\rightarrow \lambda^\sigma)$

\bigskip

$\displaystyle{\frac{\enspace \Gamma_1, \mathscr{A}[t^\tau], \lambda x^\tau \mathscr{A}[x^\tau](t^\tau), \Gamma_2 \enspace \rightarrow \enspace \Delta \enspace}{\enspace \Gamma_1, \lambda x^\tau \mathscr{A}[x^\tau](t^\tau), \Gamma_2 \enspace \rightarrow \enspace \Delta \enspace }}$ $(\lambda^\sigma \rightarrow)$

\bigskip

The Ketonen-type inference-preserving Gentzen system for classical (resp. intuitionistic) predicate logic was introduced in Inou\'{e} in \cite{Inoue2024} (resp. in \cite{InoueK1}). 

\bigskip

In the system \(\mathsf{KCT}_{h}\), the notions of a \emph{proof figure} (derivation) and of \emph{provability} are the same as in \(\mathsf{KCT}\).

\section{Ketonen-type inference-preserving tableau system for classical simple type theory $\bf KCTT_h$}

Here we introduce a Ketonen-type inference-preserving tableau system for classical simple type theory $\bf KCTT_h$ as follows.

\bigskip

 
Here we shall introduce a Ketonen-type Gentzen classical simple type theory $\bf KCTT_h$ as follows.
\medskip

Axioms: $\Gamma_1, P, \Gamma_2 \enspace \rightarrow \enspace \Delta_1, P, \Delta_2,$
\medskip

\noindent where $P$ is atomic and all A-formulas and S-formula of $\Gamma_1, \Gamma_2, \Delta_1, \Delta_2$ are atomic and $\Gamma_1, \Gamma_2, \Delta_1, \Delta_2$ may be empty, respectively.
\medskip

Rules of inference:
\bigskip

$\displaystyle{\frac{\enspace \Gamma_1, \Gamma_2 \enspace \rightarrow \enspace \Delta_1, A \supset B, \Delta_2}{\enspace \Gamma_1, A, \Gamma_2 \enspace \rightarrow \enspace \Delta_1, B, A \supset B, \Delta_2 \enspace}}$ $ (\rightarrow \supset)$.

\bigskip

$\displaystyle{\frac{\enspace \Gamma_1, A \supset B, \Gamma_2 \enspace \rightarrow \enspace \Delta \enspace }{\Gamma_1, A \supset \forall x^1x^1, A \supset B, \Gamma_2 \rightarrow \Delta \enspace \enspace \Gamma_1, B, A \supset B, \Gamma_2 \rightarrow \Delta \enspace }}$ $ (\supset \rightarrow)$,

\medskip

\noindent if $B$ is not a formula $\forall x^1x^1$.

\bigskip

$\displaystyle{\frac{\enspace \Gamma \enspace \rightarrow \enspace \Delta_1, \forall x^\sigma \mathscr{F}[x^\sigma], \Delta_2 \enspace }{\enspace \Gamma \enspace \rightarrow \enspace \Delta_1, \mathscr{F}[a^\sigma],  \forall x^\sigma \mathscr{F}[x^\sigma], \Delta_2 \enspace}}$ $ (\rightarrow \forall^\sigma)$,

\bigskip

\noindent provided that the variable $a^\sigma$ does not occur in the lower sequent.

\bigskip

$\displaystyle{\frac{\enspace \Gamma_1,\forall x^\sigma \mathscr{F}[x^\sigma], \Gamma_2 \enspace \rightarrow \enspace \Delta \enspace}{\enspace \Gamma_1, \mathscr{F}[t^\sigma], \forall x^\sigma \mathscr{F}[x^\sigma], \Gamma_2 \enspace \rightarrow \enspace \Delta \enspace}}$ $ (\forall^\sigma \rightarrow)$,
\bigskip

\noindent where $t^\sigma$ is a term.

\bigskip

$\displaystyle{\frac{\enspace \Gamma \enspace \rightarrow \enspace \Delta_1, \lambda x^\tau \mathscr{A}[x^\tau](t^\tau), \Delta_2 \enspace}{\enspace \Gamma \enspace \rightarrow \enspace \Delta_1, \mathscr{A}[t^\tau], \lambda x^\tau \mathscr{A}[x^\tau](t^\tau), \Delta_2 \enspace }}$ $(\rightarrow \lambda^\sigma)$

\bigskip

$\displaystyle{\frac{\enspace \Gamma_1, \lambda x^\tau \mathscr{A}[x^\tau](t^\tau), \Gamma_2 \enspace \rightarrow \enspace \Delta \enspace }{\enspace \Gamma_1, \mathscr{A}[t^\tau], \lambda x^\tau \mathscr{A}[x^\tau](t^\tau), \Gamma_2 \enspace \rightarrow \enspace \Delta \enspace}}$ $(\lambda^\sigma \rightarrow)$

\bigskip

The Ketonen-type inference-preserving tableau system for classical (resp. intuitionistic) predicate logic was introduced in Inou\'{e} in \cite{Inoue2024} (resp. in \cite{InoueK1}). 

\bigskip

In the system \(\mathsf{KCTT}_{h}\), the notions of a \emph{tableau},
and of a branch or tableau being \emph{closed}, \emph{open}, or \emph{complete}
(in the tableau sense) are defined as in \(\mathsf{KCTT}\).

\section{Reduction Chains and Partial Valuations for Sequents}
\subsection{Reducible, Critical, Primitive Sequents}

We first assume that for each type $\sigma$ an enumeration $\alpha_0^{\sigma},  \alpha_1^{\sigma},  \alpha_2^{\sigma},  , \dots$ of all free variables of type $\sigma$ and an enumeration 
$s_0^{\sigma},  s_1^{\sigma},  s_2^{\sigma},  , \dots$ of all terms of type $\sigma$ are given.
\begin{definition}
The \textbf{reducible formulas} of a sequent $S$ are defined to be one of the following \textup{(1)-(3)}:

\textup{(1)} $A$ is the S-formulas of $S$ which are of the form $\forall x^{\sigma} \mathcal{F}[x^{\sigma}]$ or $\lambda x^{\tau} \mathcal{A} [x^{\tau}](t^{\tau})$, 

\textup{(2)} $B$ is the A-formulas of $S$ which are of the form $\lambda x^{\tau} \mathcal{A} [x^{\tau}](t^{\tau})$,

\textup{(3)} $A \supset B$ is S-formula of $S$, where $B$ is not a formula $\forall x^1 x^1$. 

A sequent is said to be \textbf{reducible} \rm if it has at least one reducible formula. 
\end{definition}
\begin{definition} By the \textbf{distinguished formula} of a reducible sequent $S$ we mean the reducible formula of $S$ which occurs furthest to the right in $S$.
\end{definition}
\begin{definition}
By the \textbf{critical formula} of a sequent $S$ we mean the A-formula of $S$ which are of the form $\forall x^{\sigma} \mathcal{F}[x^{\sigma}]$. A sequent is said to be \textbf{critical} if it is not reducible and has at least one critical formula.
\end{definition}

\begin{definition}
A sequent is said to be \textbf{primitive} if it is neither reducible nor critical. This is obviously the case if, and only if, every minimal A-formula and every minimal S-formula of the sequent is an atomic formula.
\end{definition}

\subsection{Definition of Reduction Chains for Sequents}
\medskip
\begin{definition}
By a \textbf{R-chain (reduction chain)} of a sequent $S$ we mean a sequence of sequents
$$S_0, S_1, S_2, \ldots$$
formed as follows:
\begin{itemize}
  \item[R1.] The initial sequent $S_0$ of the R-chain is the sequent $S$.
  \item[R2.] If a formula $S_m$ of the R-chain is an axiom or an atomic sequent, then it is the last sequent of the R-chain. We then say that the R-chain has length $m$.
  \item[R3.] If a sequent $S_m$ of the R-chain is a reducible or critical sequent but not an axiom, then $S_m$ has an immediate successor $S_{m+1}$ in the R-chain, where $S_{m+1}$ is determined by $S_m$ as follows:
  \begin{itemize}
    \item[R3.1.] If $S_m$ is a reducible sequent $\Gamma_1, A \supset B, \Gamma_2 \rightarrow \Delta$ with distinguished A-formula $(A \supset B)$, then $S_{m+1}$ is either the sequent $\Gamma_1, A \supset \forall x^1 x^1, \Gamma_2 \rightarrow \Delta$ or the sequent $\Gamma_1, B, \Gamma_2 \rightarrow \Delta$.
    \item[R3.2.] If $S_m$ is a reducible sequent $\Gamma \rightarrow \Delta_1, \forall x^{\sigma} \mathcal{F}[x^{\sigma}], \Delta_2$ with distinguished S-formula $\forall x^{\sigma} \mathcal{F}[x^{\sigma}]$, then $S_{m+1}$ is the sequent $\Gamma \rightarrow \Delta_1, \mathcal{F} \alpha^{\sigma}_i, \Delta_2$ where $i$ is the least number such that $\alpha^{\sigma}_i$ does not occur in $S_m$.
    \item[R3.3.1.] If $S_m$ is a reducible sequent $\Gamma_1, \lambda x^{\tau} \mathcal{A}[x^{\tau}] (t^{\tau}), \Gamma_2 \rightarrow \Delta$ with distinguished A-formula $\lambda x^{\tau} \mathcal{A}[x^{\tau}] (t^{\tau})$, then $S_{m+1}$ is the sequent $\Gamma_1, \mathcal{A}[t^{\tau}], \Gamma_2 \rightarrow \Delta$. 
        \item[R3.3.2] If $S_m$ is a reducible sequent $\Gamma \rightarrow \Delta_1, \lambda x^{\tau} \mathcal{A}[x^{\tau}] (t^{\tau}), \Delta_2$ with distinguished S-formula $\lambda x^{\tau} \mathcal{A}[x^{\tau}] (t^{\tau})$, then $S_{m+1}$ is the sequent $\Gamma \rightarrow \Delta_1, \mathcal{A}[t^{\tau}], \Delta_2$.
    \item[R3.4.] If $S_m$ is a critical sequent 
     \[
    \Gamma_1, \forall x^{\sigma_1}_1 \mathcal{F}_1 [x^{\sigma_1}_1], \ldots, \forall x^{\sigma_k}_1 \mathcal{F}_k [x^{\sigma_k}_k], \Gamma_2 \rightarrow \Delta_1, A, \Delta_2
    \]      
    with critical formulas
         \[
    \forall x^{\sigma_1}_1 \mathcal{F}_1 [x^{\sigma_1}_1], \ldots, \forall x^{\sigma_k}_1 \mathcal{F}_k [x^{\sigma_k}_k],
    \] 
    then $S_{m+1}$ is the sequent
    \[\Gamma_1, \Gamma_2 \rightarrow \Delta_{Suc} , \] where $\Delta_{Suc} = $ 
    \[\Delta_1, 
   (A_1 \supset (\ldots \supset (A_{(m+1) \cdot k} \supset A) \ldots)), \Delta_2
    \]
    with  $A_{i+j.k} : = \mathcal{F}_i [s^{\sigma_i}_j]$
    $(i = 1, \dots ,k ; j = 0, \dots, m)$.
  \end{itemize}
\end{itemize}
\end{definition}

\bf Comment. \rm Our reduction chains (R-chains) is nothing but Sch\"{u}tte's D-chains. Sch\"{u}tte's terminology seems to be unsuitable if one thinks of its substantial meaning. Instead of D-chains we shall use reduction chains (R-chains) as the meaningly correct usage.
\medskip

\subsection{Partial Valuations for Formulas}
\medskip
\begin{definition}
By a \textbf{partial valuation (for formulas)} we mean a function $V$ which assigns a truth value $VF_i = t$ or $VF_i = f$ to certain formulas $F_i$ (but which may not be defined for all formulas) such that the following conditions are satisfied:

\begin{itemize}
  \item[V1.] If $V(A \supset B) = t$, then $VA = f$ or $VB = t$. (In this case, $VB$ or $VA$ may be undefined.)
  \item[V2.] If $V(A \supset B) = f$, then $VA = t$ and $VB = f$.
  \item[V3.] If $V^{x^\sigma} \mathcal{F}[x^\sigma] = t$, then $V \mathcal{F}[t^\sigma] = t$ for every term $t^\sigma$ of type $\sigma$.
  \item[V4.] If $V^{x^\sigma}\mathcal{F}[x^\sigma] = f$, then there is a free variable $\alpha^\sigma$ of type $\sigma$ such that $V \mathcal{F}[\alpha^\sigma] = f$. (In this case, $V \mathcal{F}[t^\sigma]$ need not be defined for every term $t^{\sigma}$ of type $\sigma$.)
    \item[V5.] If $V \lambda x^\tau \mathscr{A} [x^\tau] (t^\tau)$ is defined then $V \mathscr{A} [t^\tau] = V \lambda x^\tau \mathscr{A} [x^\tau] (t^\tau)$.
\end{itemize}
\end{definition}

\noindent
\textbf{Remark.} These partial valuations correspond to the semivaluations in Sch\"{u}tte {schutte1951}, Takahashi \cite{takahashi1967}, and Prawitz \cite{prawitz1968}.
\medskip

\subsection{Partial Valuations for Sequents}
\medskip
\begin{definition}
By a \textbf{partial valuation (for sequents)} we mean a function $V^s$ which assigns a sequent $V \Gamma \rightarrow V \Delta$ $(= V^s(S)$, $S = \Gamma \rightarrow \Delta)$, 
where $V\Gamma = VA_1, \dots , VA_n$ for $\Gamma = A_1, \dots , A_n (n \geq 0) $ and $V\Delta = VB_1, \dots , VB_m$ for $\Delta = B_1, \dots , B_m (m \geq 0)$.
\end{definition} 
\begin{definition}
By \textbf{the truth value $tv(S)$ of sequent} $S$ of truth values we mean the truth value calculated as the corresponding formula of $S$.
\end{definition}

\bf Examples. \rm 

(1) $S = t, f, t \rightarrow f, t, f$, 

$\qquad tv(S) = (t \wedge f \wedge t) \supset (t \vee f \vee t) = f  \supset t = t$.
\smallskip

(2) $S = t, f, t, f \rightarrow $,

$\qquad tv(S) =  \neg (t \wedge f \wedge t \wedge f) = \neg  f = t$.
\smallskip

(3) $S = \enspace \rightarrow f, f, f$,

$\qquad tv(S) =  f \vee f \vee f = f$.
\smallskip

(4) $S = \enspace \rightarrow $,

$\qquad tv(S) =  f$.
\smallskip

\section{Principal Lemmata for $\bf KCTT_h$}
In this section, we present a number of Lemmata for $\bf KCTT_h$. However, for the other 
systems, the same theorems hold.

\begin{lemma}[Principal Syntactic Lemma for Sequents]
If every R-chain of a sequent $S$ contains an axiom then the sequent $S$ is deducible.
\end{lemma}

\noindent Proof. In this case every R-chain of $S$ is finite. Hence by König's lemma $S$ has only finitely many R-chains. These constitute a reduction of the sequent $S$. $\Box$

\begin{lemma} [Principal Semantic Lemma for Sequents]
If there is a R-chain of sequent $S$ which contains no axiom then there is a partial valuation $V^s$ such that $tv(S) = f$
\end{lemma}
Proof. Suppose
$$(\text{RC}) \quad S_0, S_1, S_2, \ldots$$
is a R-chain which contains no axiom where $S_0$ is the sequent $S$. Using the definition of Hintikka sequent, we can now prove, in analogy with §5.3 in Sch\"{u}tte's Proof Theory \cite{schutte1977}, if we understand that positive (resp. negaitve) part corresponds to S-formula (resp. A-formula).  Concretely we may put $t$ (resp. $f$) to atomic A-formula (resp. S-formula) in a given formula (sequent) in order to obtain the falsity of the given formula (sequent). $\Box$

\begin{lemma} 
If an atomic formula occurs as a S-formula (or A-formula) in a sequent $F_m$ of the R-chain, then it also occurs as a S-formula (or A-formula) in every sequent $S_m$ ($m \geq n$) of the R-chain.
\end{lemma}
Proof. This follows immediately from the definition of R-chain. $\Box$

\begin{lemma} 
If a sequent $S_n$ of the R-chain has a reducible A-formula (or reducible S-formula) $C$, then the R-chain contains a reducible sequent $S_m$, ($m \geq n$), with distinguished S-formula (or A-fromula) $C$.
\end{lemma}
Proof. By induction on the sum of ranks of reducible formula of $S_n$, which occur in $S_n$, to right side $C$. $\Box$

\begin{lemma} 
If any sequent $S_n$ from the chain has A-formula $\forall x^{\sigma} \mathcal{F}[x^{\sigma}]$ then the chain contains infinitely many critical sequent with A-formula $\forall x^{\sigma} \mathcal{F}[x^{\sigma}]$.
\end{lemma}

Proof. Under the hypothesis for Lemma 9.5, it suffices to prove that the chain contains a critical sequent $S_m$ ($m \geq n$). We call the sum of ranks the \textbf{reducibility rank} of the sequent. By hypothesis the sequent $S_n$ is either reducible or critical. The R-chain therefore contains a sequent $S_{n+1}$, which also contains the formula $\forall x^\sigma \mathcal{F} [x^\sigma]$ as an A-formula. We prove our assertion by induction on the \textbf{reducibility rank} of the sequent $S_{n+1}$. 

If $S_{n+1}$ is a critical sequent then the assertion is satisfied with $m = n+1$. Otherwise, the R-chain contains a sequent $S_{n+2}$ which has a smaller reducibility rank than $S_{n+1}$ and also contains the formula $\forall x^\sigma \mathcal{F}[x^\sigma]$ as an A-fromula. In this case, the assertion follows from the induction hypothesis. $\Box$

Now we set $VC = f$ (or $VC = f$) if $C$ occurs in the R-chain (RC) as a S-formula (or S-formula).

\begin{lemma}
The function $V$ satisfies conditions V1–V5 for a partial valuation.
\end{lemma}

Proof.
\begin{enumerate}
  \item Suppose $V(A \supset B) = t$, then $(A \supset B)$ is an A-formula in (RC). If $B$ is a formula $\forall x^{1} x^{1}$, then $A$ is a S-formula in (RC) and therefore $VA = t$. Otherwise, it follows from Lemma 9.4 and the definition of R-chain that $A$ is a S-formula or $B$ is a A-formula in (RC), and therefore $VA = f$ or $VB = t$.

  \item Suppose $V(A \supset B) = f$, then $(A \supset B)$ is a S-formula in (RC). But then $A$ is an A-formula and $B$ a S-formula in (RC), so $VA = f$ and $VB = t$.

  \item Suppose $\forall x^\sigma \mathcal{F}[x^\sigma] = t$, then $\forall x^\sigma \mathcal{F}[x^\sigma]$ is an A-formula in (RC). From Lemma 9.5 and the definition of R-chain it follows that, for every term $t^\sigma$ of type $\sigma$, $\mathcal{F}[t^\sigma]$ is an A-formula in (RC). Hence $VF[t^\sigma] = t$ for all $t^\sigma$ of type $\sigma$.

  \item Suppose $V \forall x^\sigma \mathcal{F}[x^\sigma] = f$, then $\forall x^\sigma \mathcal{F}[x^\sigma]$ is a S-formulain (RC). From Lemma 9.4 and the definition of R-chain it follows that there is a variable $a^\sigma$ of type $\sigma$ such that $\mathcal{F}[a^\sigma]$ is a pS-formula in (RC), and therefore $V\mathcal{F}[a^\sigma] = f$.

  \item Suppose $V \lambda x^\tau \mathcal{A} [x^\tau][t^\tau]$ is defined. Then $\lambda x^\tau \mathcal{A} [x^\tau][t^\tau]$ is a S-formula (or A-formula) in (RC). From Lemma 9.4 and the definition of R-chain it follows that $\mathcal{A}[t^\tau]$ is a S-formula (or A-formula)t in (RC). Hence
  \[
  V\mathcal{A}[t^\tau] = V\lambda x^\tau \mathcal{A}[x^\tau][t^\tau].
  \]
\end{enumerate}

Since the R-chain (RC) contains no axiom, it follows from Lemma 9.3 that no atomic formula occurs both as a S-formula and as an A-formula in (RC). Therefore, there is no atomic formula $P$ for which $VP = t$ and $VP = f$.
From Lemma 9.6, it follows by induction on the rank of the formula $C$ that there is no formula $C$ for which $VC = t$ and $VC = f$. 

Hence $V$ is a partial valuation. 

Hence $VF = f$, since $F$ occurs as a S-formula in (RC). $\Box$

\begin{lemma}
The function $V^s$ for sequent satisfies conditions V1–V5 for a partial valuation.
\end{lemma}

\noindent Proof. Immediate form Definition 8.6 and definition 8.7 and Lemma 9.6. $\Box$

\section{Hintikka sequents for $\bf KCTT_h$}

\begin{definition}
Let $\Gamma$ be a sequence of formulas. We define a \it sequence set \rm $\Gamma_s$ of $\Gamma$ as follows:
$$\Gamma_s := \{ A \mid A \text{ occurs in } \Gamma \}$$
\end{definition}


\begin{definition}[Sequent set pair]
Let $S$ be a sequent, say $S := \Gamma \rightarrow \Delta$. A \emph{sequent set pair} $S_s$ is an ordered pair $(\Gamma_s, \Delta_s).$
\end{definition}

\begin{definition}[Saturation of sequence set pair]
Let $S$ be a sequent, say $S := \Gamma \rightarrow \Delta$. 
A sequence set pair $(\Gamma_s,\Delta_s)$ of $S$ is said to be \textbf{saturated}
if it satisfies the following conditions.

\begin{enumerate}
\item[(H1)]
$A \supset B \in \Delta_s \Rightarrow (A \in \Gamma_s \ \text{and}\ B \in \Delta_s)$.

\item[(H2)]
$(A \supset B \in \Gamma_s \text{ and $(B$ is not the formula $\forall x^1 x^1$}))  \Rightarrow (A \supset \forall x^1 x^1 \in \Gamma_s \ \text{or}\ B \in \Gamma_s).$

\item[(H3)]
$\forall x^\sigma\, \mathscr{F}[x^\sigma] \in \Delta_s \Rightarrow (\mathscr{F}[a^\sigma] \in \Delta_s \text{ and $a^\sigma$ occur  only in the form of formula $\mathscr{F}[a^\sigma]$}).$

\item[(H4)]
$\forall x^\sigma\, \mathscr{F}[x^\sigma] \in \Gamma_s \Rightarrow  (\text{there exists a term $t^\tau$ such that } \mathscr{F}[t^\sigma] \in \Gamma_s)$.

\item[(H5)]
$\lambda x^\tau\, \mathscr{A}[x^\tau](t^\tau) \in \Delta_s \Rightarrow \mathscr{A}[t^\tau] \in \Delta_s$.

\item[(H6)]
$\lambda x^\tau\, \mathscr{A}[x^\tau](t^\tau) \in \Gamma_s \Rightarrow \mathscr{A}[t^\tau] \in \Gamma_s$.

\end{enumerate}
\medskip

\noindent We call the (H1)--(S6) to be \textbf{saturation properties}.

\end{definition}

\textbf{Remark.} The saturation properties correspond to the reduction rules of $\bf KCTT_h$. 

\begin{definition}[Hintikka sequent]
Let $S$ be a sequent, say $S := \Gamma \rightarrow \Delta$. 
A sequent $S$ is a \textbf{Hintikka sequent} if the following conditions are satisfied:

\textup{(1)} There is a sequent $\rightarrow A$ such that $S$ the last element of R-chain of the sequent $\rightarrow A$.

\textup{(2)} The sequent set pair of $S$ is saturated in each step of reductions corresponding saturation property.. 

\textup{(3)} There is no formula $A$ such that $A \in \Gamma_s$ and $A \in \Delta_s$.

\textup{(4)} Let the first element of R-chanin of $S$ be $\rightarrow A$. Let the sequence set pair of $S$  be $(\Gamma_s,\Delta_s)$. Then for every atomic subformula $B$, $B \in \Gamma_s \cup \Delta_s$ holds.
\end{definition}

\begin{lemma}[Hintikka Sequent Lemma]
Let $\mathcal{T}$ be a completed tableau for $\bf KCTT_h$. If the tableau is open, then there is an open branch whose bottom, say $S := \Gamma \rightarrow \Delta$, is a Hintikka sequent.
\end{lemma}

\noindent Proof. It is straightfoward from the reduction rules of $\bf KCTT_h$, Definitions 10.1--10.4  and the definition of R-chain. $\Box$
\medskip

\section{Deducible sequents and Permissible Inferences for $\bf KCTT_h$}
The following theorems in this section are proved in a standard way to prove. So we omit the proofs.

\begin{theorem}[Axiom Theorem] Every sequent of the form $\Gamma_1, A, \Gamma_2 \enspace \rightarrow \enspace \Delta_1, A, \Delta_2$ is deducible.
\end{theorem}

\begin{theorem} Every sequent $\Gamma_1, \forall x^{1} x^{1}, \Gamma_2 \enspace \rightarrow \enspace \Delta$ is deducible.
\end{theorem} 

\begin{theorem} [Inversion rules] The following are weak inferences.

\textup{(1)} $(\vdash_{\bf KCTT_h} \enspace \Gamma_1, A \supset \forall x^{1} x^{1}, A \supset B, \Gamma_2 \rightarrow \Delta \enspace$ $ \mbox{and} \enspace$ $\vdash_{\bf KCTT_h}  \enspace \Gamma_1, B, \Gamma_2 \rightarrow \Delta)$ 

$\enspace \mbox{implies}$ $\vdash_{\bf KCTT_h}  \enspace \Gamma_1, A \supset B, \Gamma_2 \rightarrow \Delta.$

\textup{(2)} $\vdash_{\bf KCTT_h} \enspace \Gamma \rightarrow \Delta_1, \mathscr{F}[a^{\sigma}], \forall x^{\sigma} \mathscr{F}[x^{\sigma}], \Delta_2$ $ \mbox{implies}$
$\vdash_{\bf KCTT_h}  \enspace \Gamma \rightarrow \Delta_1, \forall x^{\sigma} \mathscr{F}[x^{\sigma}] , \Delta_2.$

\textup{(3)} $\vdash_{\bf KCTT_h} \enspace \Gamma_1, \mathscr{A}[t^\tau], \Gamma_2 \rightarrow \Delta$ $\enspace \mbox{implies}$
$\vdash_{\bf KCTT_h}  \enspace \Gamma_1, \lambda x^{\tau}\mathscr{A}[x^{\tau}(t^\tau)], \Gamma_2 \rightarrow \Delta.$

\textup{(4)} $\vdash_{\bf KCTT_h} \enspace \Gamma \rightarrow \enspace \Delta_1, \mathscr{A}[t^\tau], \Delta_2$ $\enspace \mbox{implies}$
$\vdash_{\bf KCTT_h}  \enspace \Gamma \rightarrow \Delta_1, \lambda x^{\tau}\mathscr{A}[x^{\tau}(t^\tau)], \Delta_2.$
\end{theorem} 

$S \vdash^s T$ (“S follows structurally from T”) denotes that every minimal S-formula of $S$ also occurs as a S-formula of $T$ and every minimal A-formula of $S$ also occurs as a A-formula of $T$.

\begin{theorem} [Structural rule of inference] If $F \vdash^{S} G$ holds, then $F \vdash G$ is a weak inference.
\end{theorem} 

\begin{theorem} The following is a permissible inference:

$\vdash_{\bf KCTT_h} \enspace \Gamma_1, \mathscr{F}[t^{\sigma}], \Gamma_2 \rightarrow \enspace \Delta$ $\enspace \mbox{implies}$
$\vdash_{\bf KCTT_h}  \enspace \Gamma_1, \forall x^{\sigma} \mathscr{F}[x^{\sigma}], \Gamma_2 \rightarrow \Delta.$
\end{theorem}

\begin{theorem} [General substitution rule] $ $

\textup{(1)} If $a^{\sigma}$ does not occur in the nominal form $\mathscr{F}$, then the following is permissible inference:

$\vdash_{\bf KCTT_h} \enspace \Gamma_1, \mathscr{F}[a^{\sigma}], \Gamma_2 \rightarrow \enspace \Delta$ $\enspace \mbox{implies}$
$\vdash_{\bf KCTT_h}  \enspace \Gamma_1, \forall x^{\sigma} \mathscr{F}[t^{\sigma}], \Gamma_2 \rightarrow \Delta.$

\textup{(2)} If $a^{\sigma}$ does not occur in the nominal form $\mathscr{F}$, then the following is permissible inference:

$\vdash_{\bf KCTT_h} \enspace \Gamma \rightarrow \Delta_1, \mathscr{F}[a^{\sigma}], \Delta_2$ $\enspace \mbox{implies}$
$\vdash_{\bf KCTT_h}  \enspace \Gamma \rightarrow \Delta_1, \mathscr{F}[t^{\sigma}] , \Delta_2.$
\end{theorem} 

\section{Semantics for $\bf KCTT_h$}

\subsection{Total Valuations over a System of Sets}
We consider our semantics for $\bf KCTT_h$.

\begin{definition}
By \textbf{a system $M$ of sets} we mean a collection of non-empty sets $M^{\sigma}$ for each type $\sigma$, such that $M^{\sigma_1}$ and $M^{\sigma_2}$ are disjoint if $\sigma_1 \neq \sigma_2$.
\end{definition}

For a given system $M$ of sets we introduce the following terminology.

\begin{definition} $ $

\textup{(1)} The \textbf{basic $M$-terms} of type $\sigma$ are the names of the elements of $M^{\sigma}$.

\textup{(2)} The \textbf{$M$-terms} of type $\sigma$ are the expressions which result from terms of type $\sigma$ when all the free variables occurring in them are replaced by basic $M$-terms of the same type.

\textup{(3)}The \textbf{$M$-formulas} are the \textbf{$M$-terms} of type $1$.

\textup{(4)} By an \textbf{$M$-variant} of a formula $F$ we mean an $M$-formula which is obtained from $F$ by replacing the free variables occurring in $F$ by $M$-terms of the same type where a free variable which occurs more than once in $F$ is replaced at all occurrences by the same $M$-term.

\textup{(5)}The \textbf{$M$-sequent} is similarly understood for a sequent as (3)

\textup{(6)}The \textbf{$M$-variant} of a sequent is similarly understood for a sequent as (4).
\end{definition}

As \textbf{syntactic symbols} we use\\  
$c^{\sigma}, c^{\tau_i}_{i}$ \quad for basic $M$-terms of type $\sigma, \tau_i$,\\  
$c^\tau$ \quad for sequences $c^{\tau_1}, ..., c^{\tau_n}$ $(n \geq 1)$ of basic $M$-terms,\\  
$u^\sigma, u^{\tau_i}_i$ \quad for $M$-terms of type $\sigma, \tau_i$,\\  
$u^\tau$ \quad for sequences $u^{\tau_1}_i, \dots , u^{\tau_n}_n$ $(n \geq 1)$ of $M$-terms,\\  
$A', B', C', F', G', \forall x^\sigma \mathscr{F}' [x^{\sigma}]$ \quad for $M$-formulas and\\  
$\lambda x^\tau \mathscr{A}' [x^{\tau}]$ \quad for $M$-terms of type $(\tau) = (\tau_, \dots , \tau_n)$.

\begin{definition}
By a \textbf{total valuation} over a system $M$ of sets we mean a function $W$ which assigns a truth value $WF' = t$ or $WF' = f$ to each $M$-formula $F'$ such that the following hold:

W1. \quad $W(A' \supset B') = t$ if, and only if $WA' = f$ or $WB' = t$.

W2. \quad $W \forall x^\sigma {\mathscr{F}}' [x^{\sigma}] = t$ if, and only if $W {\mathscr{F}}' [u^\sigma] = t$ for every $M$-term $u^\sigma$ of type $\sigma$.

W3. \quad $W \lambda x^{\tau} {\mathscr{A}}' [x^{\tau}] (u^\tau) = W {\mathscr{A}}' [u^{\tau}]$.

A formula $F$ of $\bf KCTT_h$ is said to be true in a total valuation $W$ over $M$ if $WF' = t$ for every $M$-variant $F'$ of $F$.
\end{definition}

\begin{definition}
Let $W$ be a \textbf{total valuation} over a system $M$ of sets. 
By a \textbf{total valuation (for sequents)} we mean a function $W^s$ which assigns a sequent $W \Gamma \rightarrow W \Delta$ ($= W^s(S)$, $S = \Gamma \rightarrow \Delta$), 
where $W\Gamma = WA_1, \dots , WA_n$ for $\Gamma = A_1, \dots , A_n (n \geq 0) $ and $W\Delta = WB_1, \dots , WB_m$ for $\Delta = B_1, \dots , B_m (m \geq 0)$.
\end{definition} 

\begin{definition}
A formula of $\bf KCTT_h$ is said to be \textbf{valid} if it is true in every total valuation over a system of sets.
\end{definition}

We follow the following convention.
\medskip

\textbf{Truth Value Convention.} We also use the truth value function $tv$ for total valuation with respect to sequents (cf. Definition 8.8).

\subsection{Soundness Theorem}

\smallskip

We here suppose every sequent $S$ of $\bf KCTT_h$ is satisfied with the following convention.
\medskip

\textbf{Name Convention}: We denote by $S'$ when the free variables in $S$ are replaced by basic $M$-terms of the same types.
\medskip

For every total valuation $W$ over $M$ we obviously have:

\begin{lemma}[SA-Lemma] 
Let $S := \Gamma_1, A, \Gamma_2 \rightarrow \Delta$ and $T := \Gamma \rightarrow\Delta_1, A, \Delta_2$.

\textup{(1)} If $WA' = f$, then $tv(A') = t$ for $S'$.

\textup{(2)} If $WA' = t$, then $tv(A') = t$ for $T'$.

\end{lemma}

\begin{theorem}[Soundness Theorem] Every deducible formula in $\bf KCTT_h$ is valid.
\end{theorem}
\noindent Proof. Let $S$ be a sequent such that $\rightarrow F$.
 Let $F$ be a deducible formula, $W$ a total valuation over a system $M$ of sets and $F'$ an $M$-variant of $F$. We prove that $WF' = t$ by induction on the order of a reduction of $F$.
 
 \quad
1. Suppose $S$ is an axiom $\Gamma_1, P, \Gamma_2 \rightarrow \Delta_1, P, \Delta_2$. Then $S'$ is an $M$-sequent $\Gamma_1', A', \Gamma_2' \rightarrow \Delta_1', A', \Delta_2'$. By the SA-Lemma we have $WF' = t$


\quad
2. Suppose $S$ is a sequent $\Gamma_1, A \supset B, \Gamma_2 \rightarrow \Delta$, which is reduced by $(\supset \rightarrow)$-inference to $\Gamma_1, A \supset \forall x^1x^1, A \supset B, \Gamma_2 \rightarrow \Delta$ and $\Gamma_1, B, A \supset B, \Gamma_2 \rightarrow \Delta$.
Then $S'$ is an $M$-sequent $\Gamma_1', A' \supset \forall x^1x^1, A' \supset B', \Gamma_2' \rightarrow \Delta'$.
If $WA' = t$, then $W(A' \supset B') = WB'$.
Since, by I.H., $W(\Gamma_1, B, A \supset B, \Gamma_2 \rightarrow \Delta) = t$, we have $WF' = t$.
If $WA' = f$, then $W(A' \supset B') = t = W(A' \supset \forall x^1 x^{1})$.
Since, by I.H., $W(\Gamma_1', A' \supset \forall x^1x^1, A' \supset B', \Gamma_2' \rightarrow \Delta') = t$, we have $WF' = t$.

\quad
3. Suppose $S$ is a sequent $\Gamma \rightarrow \Delta_1, \forall x^\sigma \mathscr{F}[x^\sigma], \Delta_2$, which is reduced by $(\rightarrow \forall^\sigma)$-inference to $\Gamma \rightarrow \Delta_1, \mathscr{F}[a^\sigma],  \forall x^\sigma \mathscr{F}[x^\sigma], \Delta_2$. Then $S'$ is an $M$-sequent $\Gamma' \rightarrow \Delta_1', \forall x^\sigma \mathscr{F}'[x^\sigma], \Delta_2'$.
If $W \forall x^\sigma \mathscr{F}' [x^\sigma] = t$, then $WF' = t$ by the SA-Lemma.
If $W \forall x^\sigma \mathscr{F}' [x^\sigma] = f$, then there is an $M$-term $u^\sigma$ such that $W \mathscr{F}' [u^\sigma] = f$. 
 $\Gamma \rightarrow \Delta_1, \mathscr{F}[u^\sigma],  \forall x^\sigma \mathscr{F}[x^\sigma], \Delta_2$ is an $M$-variant of $\Gamma \rightarrow \Delta_1, \mathscr{F}[a^\sigma],  \forall x^\sigma \mathscr{F}[x^\sigma], \Delta_2$.
Hence, by I.H., $W (\Gamma \rightarrow \Delta_1, \mathscr{F}[u^\sigma],  \forall x^\sigma \mathscr{F}[x^\sigma], \Delta_2) = t$.
Therefore $WF' = t$.

\quad
4. Suppose $S$ is a sequent $\Gamma_1, \forall x^\sigma \mathscr{F}[x^\sigma], \Gamma_2 \rightarrow \Delta$, which is reduced by $(\forall^\sigma \rightarrow)$-inference to $\Gamma_1, \mathscr{F}[t^\sigma], \forall x^\sigma \mathscr{F}[x^\sigma], \Gamma_2 \rightarrow \Delta$.
Then $S'$ is an $M$-sequent $\Gamma_1', \forall x^\sigma \mathscr{F}'[x^\sigma], \Gamma_2' \rightarrow \Delta'$.
By I.H., we have $W ( \Gamma_1', \mathscr{F}'[u^\sigma], \forall x^\sigma \mathscr{F}'[x^\sigma], \Gamma_2' \rightarrow \Delta') = t$ for an $M$-term $u^\sigma$.
If $W \mathscr{F}' [u^\sigma] = t$, then $WF' = t$.
If $W \mathscr{F}' [u^\sigma] = f$, then $W \forall x^\sigma \mathscr{F}' [x^\sigma] = f$ too. Then $WF' = t$ by the SA-Lemma.

\quad
5. Suppose $S$ is a sequent $\Gamma \rightarrow \Delta_1, \lambda x^\tau \mathscr{A}[x^\tau](t^\tau), \Delta_2$, which is reduced by $(\rightarrow \lambda^\sigma)$-inference to $\Gamma \rightarrow \Delta_1, \mathscr{A}[t^\tau], \lambda x^\tau \mathscr{A}[x^\tau](t^\tau), \Delta_2$.
Then $S'$ is an $M$-sequent  $\Gamma \rightarrow \Delta_1, \lambda x^\tau \mathscr{A}[x^\tau](u^\tau), \Delta_2$.
By I.H., $W(\Gamma' \rightarrow \Delta_1', \mathscr{A}'[t^\tau], \lambda x^\tau \mathscr{A}'[x^\tau](t^\tau), \Delta_2') = t$.
Since $W \lambda x^\tau \mathscr{A}' [x^\tau] (u^\tau) = W \mathscr{A}' [u^\tau] = t$ it follows $WF' = t$.

\quad
6. Suppose $S$ is a sequent $\Gamma_1, \lambda x^\tau \mathscr{A}[x^\tau](t^\tau), \Gamma_2 \rightarrow \Delta$, which is reduced by $(\lambda^\sigma \rightarrow)$-inference to $\Gamma_1, \mathscr{A}[t^\tau], \lambda x^\tau \mathscr{A}[x^\tau](t^\tau), \Gamma_2 \rightarrow \Delta$.
Then $S'$ is an $M$-sequent  $\Gamma_1, \lambda x^\tau \mathscr{A}[x^\tau](u^\tau), \Gamma_2 \rightarrow \Delta$.
By I.H., $W(\Gamma_1', \mathscr{A}'[t^\tau], \lambda x^\tau \mathscr{A}'[x^\tau](t^\tau), \Gamma_2' \rightarrow \Delta') = f$.
Since $W \lambda x^\tau \mathscr{A}' [x^\tau] (u^\tau) = W \mathscr{A}'[u^\tau] = f$ it follows $WF' = t$. $\Box$

\medskip

\subsection{Extending a Partial Valuation}
We further consider our semantics for $\bf KCTT_h$.
Let $V$ be a partial valuation. We shall define a total valuation $W$ over a suitable system $M$ of sets such that no formula $F$ for which $VF = f$ is true in the total valuation $W$. We follow D. Prawitz \cite{prawitz1968}. Our presentation in this subsection 12.3 is totally indebted to Sch\"{u}tte\textup{ \cite{schutte1977}}.

\begin{definition}\textup{(\cite{schutte1977})}
\textbf{Inductive definition} of a set $Pt^\sigma$ (of \textbf{possible valuations}) for each term $t^\sigma$ (by induction on the height (see Definition 2.12) of the type $\sigma$).

\quad
P1. $Pt^{0} := \{t^{0}\}$.

\quad
P2. If $Vt^{1}$ is defined, let $Pt^{1} := \{Vt^{1}\}$.
Otherwise let $Pt^{1} := \{t, f\}$.

\quad
P3. For $\tau = \tau_{1}, \dots , \tau_{n}$ $(n \geq 1)$ $p_0 \in Pt^{(\tau)}_0$ if, and only if, the following conditions are satisfied:

\quad
P3.1. $p_0$ is a set of $n$-tuples
$$(\langle t^{\tau_1}_1, p_1 \rangle, \dots , \langle t^{\tau_n}_1, p_n \rangle)$$
such that $p_i \in Pt_i^{\tau_i}$ $(i = 1, \dots , n)$.

\quad
P3.2. If $Vt_0^{(\tau)}(t_1^{\tau_1}, \dots , t_n^{\tau_n})$ is defined and $p_i \in Pt_i^{\tau_i}$ $(i = 1, \dots , n)$, then
$$(\langle t^{\tau_1}_1, p_1 \rangle, \dots , \langle t^{\tau_n}_1, p_n \rangle) \in p_0$$
if, and only if, $Vt_0^{(\tau)}(t_1^{\tau_1}, \dots , t_n^{\tau_n}) = t$.

\smallskip

\quad
For every type $\sigma$ let $M^\sigma$ be the set of ordered pairs $\langle t^\sigma, p \rangle$ with $p \in Pt^{\sigma}$.
Obviously $M^\sigma$ is non-empty.
So the sets $M^\sigma$ form a system $M$ of sets.
We now consider this system $M$ and use syntactic symbols introduced in \S 12.1.

\quad
If $u^{\sigma}$ is an $M$-term of type $\sigma$ let $u^{\sigma *}$ be the term of type $\sigma$ obtained from $u^{\sigma}$ by replacing each basic $M$-term $\langle t_i^{\tau_i}, p_i \rangle$ by the term $t_i^{\tau_i}$.
\end{definition}

\begin{definition}\textup{(\cite{schutte1977})}
\textbf{Inductive definition} of the \textbf{degree} $du^{\sigma}$ of an $M$-term $u^{\sigma}$.\\
\quad d1. $du^{\sigma} := 0$ if $\sigma = 0$ or $u^{\sigma}$ is a basic $M$-term $u^\sigma$.\\
\quad d2. $du_0^{\tau} (u_1^{\tau_1}, \dots, u_n^{\tau_n}) := \max (du_0^{(\tau)}, du_1^{(\tau_1)}, \dots , du_n^{(\tau_n)}) + 1$.\\
\quad d3. $d(A' \supset B') := \max(dA', dB') + 1$.\\
\quad d4. $d \forall x^{\sigma}[x] \mathscr{F}' [x^\sigma] := d\mathscr{F}' [c^\sigma] + 1$, where $c^{\sigma}$ is a basic $M$-term.\\
\quad d5. $d \lambda x^\tau \mathscr{A}' [x^\tau] := d \mathscr{A} [c^\tau] + 1$, where $c^\tau$ is a sequence of basic $M$-terms.\\
\end{definition}

\begin{definition}\textup{(\cite{schutte1977})}
\textbf{Inductive definition} of $Wu^\sigma$ for each $M$-term $u^{\sigma}$ (by induction on the degree $du^{\sigma}$).\\
\quad W1. $Wu^0 := u^{0*}$.\\
\quad W2. $W \langle t^\sigma, p\rangle := p$ for $p \in Pt^{\sigma}$. (If $\sigma = 0$, we have $p = t^0$, hence $W \langle t^0, p \rangle = t^0 = \langle t^0, p \rangle *$ by W2, in accordance with W1.)\\
\quad W3. $Wu^{(\tau)}_0 (u_{1}^{\tau_{i}}, \dots , u_{n}^{\tau_{n}}) := t$, if\\
$$(\langle u_{1}^{\tau_1*}, Wu_{1}^{\tau_1} \rangle, \dots , \langle u_{n}^{\tau_n*}, Wu_{n}^{\tau_n} \rangle) \in Wu^{(\tau)}_0.$$
Otherwise let $Wu_0^{(\tau)} (u_{1}^{\tau_{i}}, \dots , u_{n}^{\tau_{n}}) = f$.\\
\quad W4. $W(A' \supset B') := t$, if $WA' = f$ or $WB' = t$.
Otherwise let $W(A' \supset B') := f$.\\
\quad W5. $W \forall x^{\sigma} \mathscr{F}' [x^\sigma] := t$, if $W \mathscr{F}' [c^\sigma] = t$ holds for every basic $M$-term $c^\sigma$ of type $\sigma$.
Otherwise let $W \forall x^{\sigma} \mathscr{F}' [x^\sigma] := f$.\\
\quad W6. For $\tau = \tau_1, \dots , \tau_n (n \geq 1)$ let $W \lambda x^\tau \mathscr{A}' [x^\tau]$ be the set of $n$-tuples ($c_1^{\tau_1} \dots c_n^{\tau_n}$) satisfying the condition $W \mathscr{A} [c_1^{\tau_1} \dots c_n^{\tau_n}] = t$.
\end{definition}

\begin{lemma}\textup{(\cite{schutte1977})} $Wu^{\sigma} \in Pu^{\sigma *}$.
\end{lemma}

\noindent Proof. Proof by induction on the degree $du^{\sigma}$.\\
\quad
1. Suppose $\sigma = 0$.
Then we have $Wu^0 = u^{0*}$ and $Pu^{0*} = \{ u^{0*} \}$, so $Wu^{0*} \in Pu^{0*}$.\\
\quad
2. Suppose $u^{\sigma}$ is a basic $M$-term $ \langle t^{\sigma}, p \rangle$.
Then $Wu^{\sigma} = p \in P t^{\sigma}$ and $u^{\sigma *} = t^{\sigma}$, so $Wu^{\sigma} \in Pu^{\sigma *}$.\\
\quad
Suppose $\sigma = 1$ and $u^1$ is not a basic $M$-term.
By definition $W u^1 \in \{t, f\}$.
We now only need to prove that: If $V u^{1*}$ is defined, then $W u^1 = V u^{1*}$.\\
\quad
3.1 Suppose $u^1$ is an $M$-formula $u^{(\tau)}_0 (u_1^{\tau_1}, \dots , u_n^{\tau_n})$.
Then by I.H. $Wu_0^{(\tau)} \in Pu_0^{(\tau)*}$ and $Wu_i^{\tau_i} \in Pu_i^{\tau *} (i = 1, \dots , n)$.  
If $Vu^{1*}$ is defined then, by the definition of $Pu_0^{(\tau)*}$,  
$$(\langle u_{1}^{\tau_1*}, Wu_{1}^{\tau_1} \rangle, \dots , \langle u_{n}^{\tau_n*}, Wu_{n}^{\tau_n} \rangle) \in Wu^{(\tau)}_0.$$
if, and only if, $Vu^{1*} = t$.
Hence in this case $Wu^1 = Vu^{1*}$.\\
\quad
3.2 Suppose $u^1$ is an $M$-formula $(A' \supset B')$.
Then by I.H. $W A' \in P A'^*$ and $W B' \in P B'^*$.
If $V u^{1*} = t$, then $V A'^{*} = f$ or $V B'^{*} = t$.
Then by I.H. $W A' = f$ or $W B' = t$, and therefore $W u^1 = t$.
If $V u^{1*} = f$, then $V A'^{*} = t$ and $V B'^{*} = f$.
Then by I.H. $W A' = t$ and $W B' = f$, and therefore $W u^1 = f$.\\
\quad
3.3 Suppose $u^1$ is an $M$-formula $\forall x^{\sigma 0} \mathscr{F}' [x^{\sigma 0}]$.
Then $u^{1*}$ is a formula $\forall x^{\sigma 0} \mathscr{F} [x^{\sigma 0}]$.
By I.H. $W \mathscr{F} [c^{\sigma 0}] \in P \mathscr{F} [c^{\sigma 0*}]$ for every basic $M$-term of type $\sigma_0$.
If $V u^{1*} = t$ then $V \mathscr{F} [t^\sigma] = t$ holds for every term $t^{\sigma 0}$ of type $\sigma_0$.
Then by I.H. $W\mathscr{F}' [c^{\sigma 0}] = t$ for every basic $M$-term $c^{\sigma 0}$ of type $\sigma_0$, and therefore $Wu^1 = t$.
If $Vu^{1*} = f$, then there is a free variable $a^{\sigma 0}$ of type $\sigma_0$ such that $V \mathscr{F}[a^{\sigma 0}] = f$.
Then by I.H., $W \mathscr{F} [c^{\sigma 0}] = f$ for some basic $M$-term $c^{\sigma 0}$ of type $\sigma_0$ and therefore $Wu^1 = f$.\\
\quad
4. Suppose $u^\sigma$ is an $M$-term $\lambda x^\tau \mathscr{A}' [x^\tau]$ where $\sigma = (\tau) = (\tau_1, \dots , \tau_n)$ $(n \geq 1)$.
Then $u^{\sigma *}$ is a term $\lambda x^\tau \mathscr{A} [x^\tau]$.
If $Vu^{\sigma *} (c_1^{\tau_1 *}, \dots , c_n^{\tau_n *})$ is defined then $V \mathscr{A} [c_1^{\tau_1 *}, \dots , c_n^{\tau_n *}] = Vu^{\sigma *} (c_1^{\tau_1 *}, \dots , c_n^{\tau_n *})$.
Then by I.H., $W \mathscr{A} [c_1^{\tau_1}, \dots , c_n^{\tau_n}] = V \mathscr{A} (c_1^{\tau_1 *}, \dots , c_n^{\tau_n *})$.
In this case, by the definition of $Wu^\sigma$,
$$(c_1^{\tau_1}, \dots , c_n^{\tau_n}) \in Wu^{\sigma}$$
if, and only if $Vu^{\sigma *} (c_1^{\tau_1 *}, \dots , c_n^{\tau_n *}) = t$.
Hence we have $Wu^\sigma \in Pu^{\sigma *}$ by the definition of $Pu^{\sigma *}$. $\Box$

\begin{lemma}\textup{(\cite{schutte1977})} If $\nu$ is a 1-place nominal form such that $\nu[u^\sigma]$ is an $M$-term, then
$W\nu[u^\sigma]) = W\nu[\langle u^{\sigma *}, Wu^\sigma \rangle]$.
\end{lemma}

\noindent Proof. By Lemma 12.2, $Wu^\sigma \in Pu^{\sigma *}$ and therefore $\langle u^{\sigma *}, Wu^\sigma \rangle$ is a basic $M$-term of type $\sigma$.
If $\nu = *_1$ then the assertion holds by the definition of $W \langle u^{\sigma *}, Wu^\sigma \rangle$.
In the other cases the assertion follows by induction on the length of the nominal form $\nu$. $\Box$

\begin{lemma}\textup{(\cite{schutte1977})} 
$W \forall x^{\sigma} \mathscr{F}' [x^{\sigma}] = t$ if, and only if, $W \mathscr{F}' [u^{\sigma}] = t$ for every $M$-term $u^{\sigma}$ of type $\sigma$.
\end{lemma}

\noindent Proof. If $W \mathscr{F}' [u^{\sigma}] = t$ for every $M$-term $u^{\sigma}$ of type $\sigma$ then $W \forall x^{\sigma} \mathscr{F}' [x^{\sigma}] = t$ by definition. If $W \forall x^\sigma \mathscr{F}' [x^\sigma] = t$, then $W \mathscr{F}' [\langle u^{\sigma *}, W u^{\sigma} \rangle] = t$ and by Lemma 2 $W \mathscr{F}' [u^\sigma] = t$ for every $M$-term $u^\sigma$ of type $\sigma$. $\Box$

\begin{lemma}\textup{(\cite{schutte1977})}
$W \lambda x^\tau \mathscr{A}' [x^\tau] (u^\tau) = W \mathscr{A}' [u^\tau]$.
\end{lemma}

Proof. Let $c_{i}^{\tau_i} := \langle u_{i}^{\tau_i *}, W u_{i}^{\tau_i} \rangle (i = 1, \dots , n)$ and $c^\tau := c_1^{\tau_1} , \dots , c_n^{\tau_n}$.
By definition $W \lambda x^\tau \mathscr{A}' [x^\tau] (u^\tau) = t$ if, and only if, $(c_1^{\tau_1} , \dots , c_n^{\tau_n}) \in W \lambda x^\tau \mathscr{A}' [x^\tau]$.
This is the case if, and only if, $W \mathscr{A}' [c^\tau] = t$ and then by Lemma 12.3 $W \mathscr{A}' [u^\tau] = t$. $\Box$

\begin{theorem}\textup{(\cite{schutte1977})} The restriction of $W$ to $M$-formulas is a total valuation over $M$ under which no formula $F$ is true for which $VF=f$.
\end{theorem}

\noindent Proof. By the definition of $W$ and Lemmata 12.4 and 12.5 $W$ is a total valuation over $M$.
Suppose $VF=f$. Replace each free variable $a$ occurring in $F$ by a basic M-term $\langle ^\sigma, p\rangle$ such that $p \in Pa^\sigma$ to obtain an $M$-formula $F'$.
Then $F'$ is an $M$-variant of $F$ with $F'^{*}=f$.
By Lemma 12.2, $WF' \in PF$.
Since $VF = f$ we have $WF' = f$.
Hence $F$ is not true under the total valuation $W$. $\Box$

\section{Completeness Theorem and Cut Rule for $\bf KCTT_h$}

\begin{theorem} [Completeness Theorem] Every valid formula is deducible.
\end{theorem}

\noindent Proof. Suppose $S$ ($\rightarrow A$) is not deducible. Then by the principal syntactic lemma (Lemma 9.1) there is a $R$-chain of $S$ which contains no axiom and therefore by the principal semantic lemma (Lemma 9.2) there is a partial valuation $V$ with $VA=f$. It follows from Theorem 12.2 that $A$ is not valid. $\Box$

The permissibility of the cut rule was first proved by W. Tait \cite{tait1965}  for second order predicate calculus (a subsystem of classical simple type theory). And M. Takahashi \cite{takahashi1967}  and D. Prawitz \cite{prawitz1968}
 proved it for the full classical simple type theory. These proofs are related to a semantic equivalent of the syntactic fundamental conjecture of Takeuti \cite{takeuti} which was developed by K. Sch\"{u}tte \cite{schutte1951}  using partial valuations.

\begin{theorem} [Cut Elimination Theorem] $ $
\medskip

$(\vdash_{\bf KCTT_h} \Gamma \rightarrow \Delta_1, A, \Delta_2$ and $ \vdash_{\bf KCTT_h} A \rightarrow B)$ implies 

$ \qquad \vdash_{\bf KCTT_h}  \Gamma \rightarrow \Delta_1, B, \Delta_2.$
\end{theorem}

\noindent Proof. Suppose that $\Gamma \rightarrow \Delta_1, A, \Delta_2$ and $\vdash A \rightarrow B$ are deducible, then they are valid by Theorem 12.1. Let $W$ be a total valuation over a system $M$ of sets.
If $\Gamma' \rightarrow \Delta_1', A', \Delta_2'$ and $\vdash A' \rightarrow B'$ are $M$-variants of them.   Then $W(\Gamma' \rightarrow \Delta_1', A', \Delta_2')=t$ and $W(A' \rightarrow B')=t$.
If $WB'=t$, then $W (\Gamma' \rightarrow \Delta_1', B', \Delta_2')=t$, that is, $tv(B')=t$ by the SA-Lemma (Lemma 12.1).
If $WB'=f$, then since $W(A' \rightarrow B')=t$, we have $WA'=f$. Then from $W(\Gamma' \rightarrow \Delta_1', A', \Delta_2')=t$ we have $W(\Gamma' \rightarrow \Delta_1', B', \Delta_2')=t$. 
In any case, therefore, $W(\Gamma' \rightarrow \Delta_1', A', \Delta_2')=t$, that is, $tv(B')=t$. Hence 
Hence $\Gamma \rightarrow \Delta_1, B, \Delta_2$ is valid and, by Theorem 13.1, deducible. $\Box$

\section{Remarks}

Tableau systems made the very way open to automatization to prove formulas. This is a great contribution of Oiva Ketonen to logic, computer science and human intellectuality (refer to Inou\'{e} and Miwa \cite{InoueMiwaKetonen}, Ketonen (eds. by Negri and von Plato) \cite{Ketonen} and Inou\'{e} \cite{Inoue2024, InoueK1, InoueK3}).

Our definition of Hintikka sequent is specific to $\bf KCTT_h$. We can present the definition in more general setting, which makes the formalization easier for  e.g. Mizar formalization.

The second author of this paper thinks that proof may be understood as sheaf. This point of view will be investigated in Inou\'{e} \cite{Inoue2026} and the study will be related to this paper in the future.

\bigskip

\noindent Tadayoshi Miwa

\noindent Library, The University of Tokyo Library

\noindent Hongo 7-3-1, Bunkyo-ku, Tokyo 113-0033, Japan

\noindent miwa.tadayoshi@mail.u-tokyo.ac.jp

\bigskip

\noindent Takao Inou\'{e}

\noindent Faculty of Informatics

\noindent Yamato University

\noindent Katayama-cho 2-5-1, Suita, Osaka, 564-0082, Japan

\noindent inoue.takao@yamato-u.ac.jp
 
\noindent (Personal) takaoapple@gmail.com (I prefer my personal mail)


\begin{thebibliography}{99}
\bibitem{InoueK4} \sc Takao  Inou\'{e}, \rm Ketonen-type Gentzen and tableau systems for classical  one-place predicate logic and its axiomatic rejection. Forthcoming.
\end{thebibliography}
\end{document}